\documentclass[12pt]{amsproc}

\theoremstyle{definition}

\theoremstyle{remark}

\numberwithin{equation}{section}

\begin{document}

\title{String Modular Motives of\\Mirrors of Rigid Calabi-Yau Varieties}

\author{Savan Kharel}
\address{Indiana University South Bend \\ South Bend, IN 46634}
\email{skharel@iusb.edu}

\author{Monika Lynker}
\address{Indiana University South Bend \\ South Bend, IN 46634}
\email{mlynker@iusb.edu}

%    Information for third author
\author{Rolf Schimmrigk}
\address{Indiana University South Bend \\ South Bend, IN 46634}
\email{netahu@yahoo.com, rschimmr@iusb.edu}

%General info
\date{July 2, 1991}

\thanks{M.L. and R.S. were supported in part by the
NSF under Grant No. PHY99-07949. They are grateful to the Kavli
Institute for Theoretical Physics, Santa Barbara, for hospitality
and support through KITP Scholarships during the course of part of
this work. R.S. was also supported through IUSB faculty research
grants. }

\begin{abstract}
 The modular properties of some higher dimensional varieties of
 special Fano type are analyzed by computing the L-function of
 their $\Omega-$motives. It is shown that the emerging modular forms
 are string theoretic in origin, derived from the characters of
 the underlying rational conformal field theory. The definition
 of the class of Fano varieties of special type is
 motivated by the goal to find candidates for a geometric realization of the
 mirrors of rigid Calabi-Yau varieties. We consider explicitly the cubic
 sevenfold and the quartic fivefold, and show that their motivic L-functions
 agree with the L-functions of their rigid mirror Calabi-Yau varieties.
 We also show that the cubic fourfold is string theoretic, with a modular form
 that is determined by that of an exactly solvable K3 surface.
\end{abstract}

\maketitle

\tableofcontents

%\parskip 0.05truein
%\baselineskip=18pt

\section{Introduction}

In this paper we continue the program of applying methods from
arithmetic geometry to the problem of understanding the question
how spacetime emerges in string theory, more precisely, the
relation between the physics on the worldsheet and the nontrivial
geometric component of spacetime. The focus here will be on a
class of varieties that have been discussed in the context of
mirror symmetry for rigid Calabi-Yau spaces  \cite{s92, cdp93,
s94}. These manifolds are higher dimensional spaces, characterized
by the existence of a nontrivial cohomology group
$H^{n-(Q-1),(Q-1)}(X)$, where $n={\rm dim}_{\mathbb C}X$ and
$Q\geq 1$ is an integer. In particular $H^{n-i,i}(X)=0$ for $0\leq
i < Q-1$ if $Q>1$, and for $Q=1$ this class of varieties reduces
to manifolds Calabi-Yau type. For $Q>1$ spaces of this type were
called Fano manifolds of special type in \cite{s92, s94}, and
generalized Calabi-Yau manifolds in \cite{cdp93}. The latter name
has in the meantime generally been adopted for Hitchin's
generalization of Calabi-Yau varieties. We will call these spaces
special Fano manifolds of charge $Q$.  In the context of exactly
solvable models manifolds of this type are distinguished by the
fact that the number of tensor factors of $N=2$ supersymmetric
minimal models exceeds the number of variables normally associated
with the standard relation between the central charge and the
dimension of the variety.

Our goal is to show that, in complete analogy to the case of
Calabi-Yau varieties considered in \cite{s06,s07} and references
therein, the arithmetic geometry of Fano varieties of special type
encodes features of the underlying conformal field theories. In
particular the L-function of the $\Omega-$motive, considered in
\cite{s06, s07b} in the context of Calabi-Yau varieties,
generalizes to the class of Fano varieties of special type. We
will show that for those members in this class that are expected
to be mirrors of rigid Calabi-Yau manifolds the inverse Mellin
transform of the L-function of the $\Omega-$motive with a Tate
twist of charge $(Q-1)$ is a modular form of weight four which
decomposes into factors derived from the underlying conformal
field theory.

In the context of finding mirrors of rigid Calabi-Yau manifolds
two varieties are of particular interest. The first is the cubic
sevenfold
 \begin{equation}
  X_7^3 = \left\{(z_0:\cdots :z_8) \in {\mathbb P}_8~{\Big |}~
    \sum_{i=0}^8 z_i^3 =0 \right\}
 \end{equation}
 whose Hodge decomposition of the intermediate cohomology is given by
 \begin{eqnarray}
   h^{7,0}(X_7^3) = h^{6,1}(X_7^3) &=& 0 \nonumber \\
         h^{5,2}(X_7^3) &=& 1 \nonumber \\
         h^{4,3}(X_7^3) &=& 84
 \label{cubic-sevenfold-hodges}\end{eqnarray}

 The second variety is the quartic fivefold
 \begin{equation}
  X_5^4 = \left\{ (z_0:\cdots :z_6) \in
  {\mathbb P}_{(1,1,1,1,1,1,2)}~{\Big |}~ \sum_{i=0}^5 z_i^4 + z_6^2
  =0\right\},
 \label{fano2}\end{equation}
 with a Hodge decomposition of the intermediate cohomology that is given by
 \begin{eqnarray}
         h^{5,0}(X_5^4) &=& 0 \nonumber \\
         h^{4,1}(X_5^4) &=& 1 \nonumber \\
         h^{3,2}(X_5^4) &=& 90.
 \label{quartic-fivefold-hodges}
 \end{eqnarray}

 The motivation for considering these hypersurfaces is that their
  cohomology suggests that they are related to conformal field theories at
  central charge $c=9$. The geometric cohomology projection of
  \cite{s92,s94} and the period computation of \cite{cdp93}
  furthermore suggest that they are of relevance in the
 context of rigid mirrors of Calabi-Yau threefolds.

 We first show that the Tate twisted $\Omega-$motives of
 these varieties are modular in terms of Hecke eigenforms of weight
 four which are
 induced by string theoretic modular forms of elliptic type. The
 latter factor into Hecke indefinite modular forms of weight one associated
 to the affine Lie algebra $A^{(1)}_1$.
 Let $X_n^d$ be a hypersurface of degree $d$ and dimension $n$,
 and denote by $L_{\Omega}^{Q-1}(X_n^d,s)$ the L-function of the $\Omega-$motives
 $M_{\Omega}(X_n^d)$ of the variety $X_n^d$ with the necessary Tate twist of
 charge $(Q-1)$. Denote by
 $\Theta^k_{\ell,m}(\tau) = \eta^3(\tau)c^k_{\ell,m}(\tau)$
 the theta functions defined in terms of the Dedekind eta function
 $\eta(\tau)$ and the Kac-Peterson string functions $c^k_{\ell,m}(\tau)$.
 The Hecke indefinite modular
 forms $\Theta^k_{\ell,m}(\tau)$ are associated to the affine Lie
 algebra $A^{(1)}_1$ of the underlying rational conformal field
 theory. Finally, we denote by $\chi_2$ denote the Legendre character,
  and set $q=e^{2\pi i\tau}$.

{\bf Theorem 1.}~
 {\it Let $(n,d) \in \{(7,3),(5,4)\}$.
 The inverse Mellin transforms $f_{\Omega}^{Q-1}(X_n^d,q)$ of the
 L-functions $L_{\Omega}^{Q-1}(X_n^d,s)$ of the $\Omega-$motive with
 a Tate twist of charge $(Q-1)$ are modular forms in
$S_4(\Gamma_0(N))$ with $N\in \{9,64\}$, described by the
algebraic Hecke characters $\psi_{d^3}$}
 \begin{equation}
  f_{\Omega}^{Q-1}(X_n^d,q) = f(\psi_{d^3}^3,q).
 \label{fanomod}\end{equation}
{\it The modular forms $f(\psi_{d^3},q) \in S_2(\Gamma_0(d^3))$
are of elliptic type, and factor into string theoretic Hecke
indefinite modular forms as}
 \begin{eqnarray}
  f(\psi_{27},q) &=& \Theta^1_{1,1}(q^3)\Theta^1_{1,1}(q^9) \nonumber \\
  f(\psi_{64},q) &=& \Theta^2_{1,1}(q^4)^2\otimes \chi_2.
 \end{eqnarray}

This result provides an interpretation of the L-function of the
twisted $\Omega-$motive of these varieties in terms of modular
forms on the string worldsheet theory of exactly solvable Gepner
models at central charge $c=9$, with nine and six minimal $N=2$
supersymmetric factors, respectively. The twist character $\chi_2$
is the quadratic character of the field of quantum dimensions of
the underlying affine Lie algebra $A_1^{(1)}$ \cite{s01, ls04,
s05}.

The characters $\psi_{27}$ and $\psi_{64}$ are the complex
multiplication Hecke characters associated to the elliptic
Brieskorn-Pham curves $E^3 \subset {\mathbb P}_2$ and $E^4 \subset
{\mathbb P}_{(1,1,2)}$ of conductor 27 and 64, respectively. These
tori are exactly solvable and their string theoretic modular forms
have been determined in \cite{su02} and \cite{ls04}, respectively,
to factor as
 \begin{eqnarray}
  f(E^3,q) &=& \Theta^1_{1,1}(q^3)\Theta^1_{1,1}(q^9) \nonumber \\
  f(E^4,q) &=& \Theta^2_{1,1}(q^4)^2\otimes \chi_2.
 \end{eqnarray}
 Their appearance in the context of the two Fano varieties $X_7^3$ and
 $X_5^4$ suggests to consider resolutions $X_3^d$
 of quotients of triple products of
elliptic curves $(E^d)^3$ for $d=3,4$, and to compare their
L-functions with those of the special Fano varieties determined in
Theorem 1. Define ${\mathbb Z}_d:= {\mathbb Z}/d{\mathbb Z}$. It
turns out that there are actions of ${\mathbb Z}_d\times {\mathbb
Z}_d$ on $(E^d)^3$ such that the resolved manifolds $X_3^d$ of the
quotient threefolds $(E^d)^3/{\mathbb Z}_d\times {\mathbb Z}_d$
have the appropriate mirror cohomology of the special Fano
varieties $X^d_{n}$ for $d=3,4$
 \begin{eqnarray}
  h^{1,1}(X_3^3) &=& 84 \nonumber \\
  h^{1,1}(X_3^4) &=& 90.
  \end{eqnarray}
 Combining Theorem 1 with results
by Cynk and Hulek \cite{ch05} leads to the following corollary.

{\bf Theorem 2.} {\it The varieties $X_3^d$, $d=3,4$, are rigid
Calabi-Yau manifolds whose L-functions of the intermediate
cohomology are modular and agree with the L-functions of the
special Fano varieties $X_n^d$ with $(d,n)\in \{(3,7),(4,5)\}$,
respectively.}

We interpret the string theoretic modularity of the special Fano
varieties $X_7^3$ and $X_5^4$, together with the agreement of
their motivic L-functions with those of their corresponding rigid
Calabi-Yau spaces, as further support for the interpretation of
certain higher dimensional Fano varieties of special type as
mirrors of rigid Calabi-Yau manifolds. It would be of interest to
extend our arithmetic analysis of the mirror pairs $(X_n^d,X_3^d)$
to families by using the methods developed in \cite{cdr00, cdr04},
and further considered in \cite{k04,k06,ky06,cdr07}.

The modular analysis described here generalizes to other
dimensions and central charges. To illustrate this we consider a
third Fano variety of special type, defined by the cubic fourfold
 \begin{equation}
  X_4^3 = \left\{(z_0:\cdots z_5) \in {\mathbb P}_5~{\Big |}~
   \sum_{i=0}^5 z_i^3=0\right\},
 \label{cubic-fourfold}\end{equation}
  whose Hodge decomposition of the intermediate cohomology is given by
 \begin{eqnarray}
   h^{4,0}(X_4^3) &=& 0 \nonumber \\
   h^{3,1}(X_4^3) &=& 1 \nonumber \\
   h^{2,2}(X_4^3) &=& 21.
 \label{cubic-fourfold-hodges}
 \end{eqnarray}
 For this variety we obtain the following string theoretic interpretation.
 Let $\vartheta(q)$ denote a modular form of weight one to be
 defined later in the paper.

 {\bf Theorem 3.}~
  {\it The inverse Mellin transform $f_{\Omega}^1(X_4^3,q)$ of the
  L-series of the Tate twisted $\Omega$-motive of the cubic fourfold
  $X_4^3$ is a modular form of
 weight three and level 27, which factors into Hecke indefinite
 modular forms as}
 \begin{equation}
  f_{\Omega}^{1}(X_4^3,q) =
  \vartheta(q^3)\Theta^1_{1,1}(q^3)\Theta^1_{1,1}(q^9)~
   \in ~S_3(\Gamma_0(27)).
 \end{equation}
 {\it This form is identical to the modular form of the $\Omega-$motive
 of the K3 surface}
  \begin{equation}
   S^6 = \left\{(z_0:\cdots z_3) \in {\mathbb P}_{(1,1,1,3)}~{\Big |}~
    z_0^6 + z_1^6 + z_2^6 + z_3^2 =0\right\},
  \end{equation}
  {\it i.e. $f_{\Omega}^1(X_4^3,q)=f_{\Omega}(S_6,q)$. It has complex multiplication
  and is determined by the Hecke character $\psi_{27}$ of
  weight 2 associated to the Eisenstein field ${\mathbb Q}(\mu_3)$}
   \begin{equation}
    f_{\Omega}^{1}(X_4^3,q)=f(\psi_{27}^2,q).
   \end{equation}
The modularity of the cubic fourfold has independently been
considered by Goto \cite{g06} and also by Hulek and Kloosterman
\cite{hk06}. The string theoretic interpretation based on the
affine algebra $A^{(1)}_1$ follows from the results of \cite{s06}.

\section{Fano varieties of special type}

In this section we briefly review the characterization of Fano
varieties of special type introduced in \cite{s92,cdp93,s94},
adopting the notation of \cite{s92,s94}. Let $Q,D_{\rm crit} \in
{\mathbb N}$ and set $s=D_{\rm crit} + 2Q -1$. We consider the
class of hypersurfaces in
 weighted projective space ${\mathbb P}_{(k_0,...,k_s)}$ with
 weights $k_i$ chosen such that
 \begin{equation}
  d = \frac{1}{Q} \sum_{i=0}^s k_i
  \end{equation}
  is an integer, defining the degree of the polynomial that
  determines the variety
 \begin{equation}
  X^d_{s-1} = \left\{(z_0:\cdots :z_s) \in {\mathbb P}_{(k_0,...,k_s)}~
     {\Big |} ~ p(z_0,...,z_s)=0, ~
     {\rm deg}~p=d \right\}.
 \label{sfano}\end{equation}
The integer $Q$ has the physical interpretation of a charge of the
underlying physical theory, and $D_{\rm crit}$ can be viewed as
the critical complex dimension of the corresponding Calabi-Yau
manifold in those cases where it can be constructed. It is
determined by the central charge $c$ of the conformal field theory
as $c=3D_{\rm crit}$.

Alternatively, varieties $X$ of this type can be characterized by
their first Chern class
 \begin{equation}
  c_1(X) = (Q-1)c_1({\mathcal N}),
 \end{equation}
 where $c_1({\mathcal N})$ is the first Chern class of the normal
 bundle ${\mathcal N}$ of $X \subset {\mathbb P}_{(k_0,...,k_s)}$.
 For $Q=1$ this construction therefore recovers the class of
 Calabi-Yau hypersurfaces for arbitrary dimensions, while for $Q>1$
 these hypersurfaces  have positive first Chern class and
 hence do not a priori describe consistent string vacua, but instead
 define a sub-class of Fano varieties. It was shown in \cite{s92,s94}
 that these varieties nevertheless are closely related to string vacua in the
 critical dimension $D_{\rm crit} = s+1-2Q$ in the sense that it is possible
 to derive the massless spectrum of critical vacua from the cohomology
 of the special Fano $(s-1)-$folds via a geometric projection.
 For particular classes of these special Fano varieties it is also
 possible to construct the critical Calabi-Yau manifold explicitly via this
 projection, in which case the charge $Q$ has the meaning of the
 codimension of the critical manifold.
 More intrinsically, we can therefore write
  \begin{equation}
   D_{\rm crit}(X) = {\rm dim}_{\mathbb C} X + 2(1-Q).
  \end{equation}

 The essential ingredient of our arithmetic analysis below is the
 existence of a nontrivial cohomology group
 $H^{n-(Q-1),(Q-1)}(X)$. We therefore use this feature as the
 characterizing property.

{\bf Definition.}~{\it A Fano variety is called special and of
charge $Q$, if it has a nonvanishing cohomology group
$H^{n-(Q-1),(Q-1)}(X)$.}

Interest in this class of manifolds arose originally because they
provide a geometric framework for $N=2$ exactly solvable Gepner
models which do not have K\"ahler deformations, hence they provide
candidates for the mirrors of rigid Calabi-Yau varieties
\cite{s92, cdp93, s94}. A toric description of these special Fano
spaces was subsequently given by Batyrev and Borisov \cite{bb97},
and Batyrev and Dais \cite{bd96}.
 A discussion of Lagrangian subvarieties in the context of the
 Strominger-Yau-Zaslow mirror conjecture \cite{syz96} was given by Leung
 \cite{l02}.

\section{String theoretic modular forms}

{\bf 3.1} The simplest class of N=2 supersymmetric exactly
solvable theories are built in terms of the affine Lie algebras.
A construction of such non-twisted affine Kac-Moody algebras, is
provided by the extension \cite{k90}
 \begin{equation}
 \hat{\underline{\rm G}} =
{\rm L}\underline{\rm G} \oplus \mathbb{C} k \oplus \mathbb{C} d
 \end{equation}
 of the loop algebra
 \begin{equation}
 {\rm L}\underline{\rm G} = \underline{\rm G} \otimes
 \mathbb{C}[t,t^{-1}]
 \end{equation}
 by the central extension $k$ and $d=t\frac{d}{dt}$. In terms of the
generators $J^a\otimes t^m$ the algebra takes the form
 \begin{equation}
 [J^a\otimes t^m,J^b\otimes t^n] = if^{ab}_cJ^c \otimes t^{m+n} +
 km \delta^{ab} \delta_{m+n,0}.
 \end{equation}
 The representations
of this algebra can be parametrized by affine roots
$\hat{\lambda}=(\lambda,k,n)$ of the Cartan-Weyl subalgebra
$\{H^i_0,E^{\alpha}_0,L_0\}$, with $i=1,...,r$, where $r$ denotes
the rank of the underlying Lie algebra $\underline{\rm G}$. For
fixed affine level $k$ of the theory the characters
$\chi_{\hat{\lambda}}$ are essentially parametrized by the weight
$\lambda$ of the representation.  For the reduced characters of an
affine Lie algebra $\underline{\hat G}$ at level $k$ the
characters transform as
 \begin{equation} \chi_{\hat{\lambda}}(-1/\tau) =
       \sum_{\hat{\mu} \in P_+^k} S_{\hat{\lambda},\hat{\mu}}
       \chi_{\hat{\mu}}(\tau),
 \end{equation}
 where the modular $S$-matrix takes
       the form
 \begin{equation} S_{\hat{\lambda},\hat{\mu}} =
        \frac{i^{|\Delta_+|}}{\sqrt{|P/Q^{\vee}| (k+g)^r}} \sum_{w\in W}
        \epsilon(w) e^{-2\pi
        i\frac{<w(\lambda+\rho),\mu+\rho>}{k+g}}.
 \end{equation}
        Here $P=\sum_i \mathbb{Z} \omega_i$ denotes the lattice
        with fundamental weights $\omega_i$ defined by
        $<\omega_i,\alpha_j^{\vee}>=\delta_{ij}$ via co-roots
        $\alpha_j^{\vee}=\frac{2\alpha}{<\alpha,\alpha>}$.
        $Q^{\vee}=\sum_i \mathbb{Z}\alpha_i^{\vee}$ is the co-root lattice
        and $P/Q^{\vee}$ denotes the lattice points of $P$ lying in an elementary
         cell of $Q^{\vee}$, while $|P/Q^{\vee}|$ describes the number of points in this
         set. $P_+^k$ is the set of all dominant weights at level
         $k$, and $\epsilon(w)=(-1)^{\ell(w)}$ is the signature of the
         Weyl group element $w$,
         where $\ell(w)$ is the minimum number of simple Weyl reflections that $w$
         decomposes into. $\Delta_{+}$ is the number of positive roots in
         the Lie algebra $\underline{G}$.

{\bf 3.2} Supersymmetric string models can be constructed in terms
of conformal field theories with $N=2$ supersymmetry. For the
present discussion the important structure is determined by the
affine Lie algebra $A_1^{(1)}$, which provides the essential
building block of the supersymmetric theory. It turns out that of
particular importance are the Hecke indefinite modular forms which
can be defined as
 \begin{equation}
  \Theta^k_{\ell,m}(\tau) =
\sum_{\stackrel{\stackrel{-|x|<y\leq |x|}{(x,y)~{\rm
or}~(\frac{1}{2}-x,\frac{1}{2}+y)}}{\in
\mathbb{Z}^2+\left(\frac{\ell+1}{2(k+2)},\frac{m}{2k}\right)}}
{\rm sign}(x) e^{2\pi i \tau((k+2)x^2-ky^2)}
\end{equation}
 Related to these forms are the string functions
 $c^k_{\ell,m}(\tau) = \Theta^k_{\ell,m}(\tau)/\eta^3(\tau)$
 of Kac-Peterson \cite{kp80,kp84}, which are of immediate physical
 relevance because they appear in the $N=2$ superconformal characters
 $\chi^k_{\ell,q,s}(\tau)$
  \begin{equation}
 \chi^k_{\ell,q,s}(\tau, z,u)
 = \sum c^k_{\ell,q+4j-s}(\tau) \theta_{2q+(4j-s)(k+2),2k(k+2)}(\tau, z,u),
 \end{equation}
 where the theta functions $\theta_m$ are defined as
 \begin{equation}
 \theta_{n,m}(\tau,z,u) = e^{-2\pi i m u} \sum_{\ell \in {\mathbb Z} +
 \frac{n}{2m}} e^{2\pi i m \ell^2 \tau + 2\pi i \ell z}.
 \end{equation}
 These characters in turn define the partition function of the conformal
 field theory.

 The theta functions $\Theta^k_{\ell,m}(\tau)$
  are associated to quadratic number fields
determined by the level of the affine theory. These are modular
forms of weight one and cannot, therefore, be identified with the
geometric modular form. It turns out, however, that appropriate
products do lead to interesting forms \cite{su02, ls04, s05, s06,
s07}.

\section{$L-$functions of $\Omega-$motives}

For a general smooth algebraic variety $X$ reduced mod $p$ the
congruence zeta function of $X/{\mathbb F}_p$ is defined by
 \begin{equation}
  Z(X/{\mathbb F}_p, t)
  = \exp\left(\sum_{r\in {\mathbb N}} \#
      \left(X/{\mathbb F}_{p^r} \right) \frac{t^r}{r}\right).
  \end{equation}
  Here the sum is over all finite extensions
  ${\mathbb F}_{p^r}$ of ${\mathbb F}_p$ of degree $r$.
 Per definition $Z(X/{\mathbb F}_p,t) \in 1 + {\mathbb Q}[[t]]$, but the
 expansion can be shown to be integer valued by writing it as an
 Euler product.

The proof by Grothendieck \cite{g65} of part of the Weil
conjectures \cite{w49} asserts that the zeta function is a
rational function determined by the cohomology of the variety
  \begin{equation}
  Z(X/{\mathbb F}_p,t) =
  \frac{\prod_{j=1}^n  {\mathcal P}_p^{2j-1}(t)}{
        \prod_{j=0}^n  {\mathcal P}_p^{2j}(t)},
   \label{ratz}\end{equation}
 where ${\rm dim}_{{\mathbb C}}X=n$, and ${\mathcal P}_p^i(t)$ is a polynomial
 ${\mathcal P}_p^i(t) = \sum_{j=0}^{b_i}\beta^i_j(p)t^j$
 associated to the $i^{\rm th}$ cohomology group with degrees $b_i$ given by the
 $i^{\rm th}$ Betti number $b_i = {\rm dim} H^i(X)$. This result motivates the
 introduction of L$-$functions associated to the polynomials
 ${\mathcal P}_p^i(t)$,
 thereby reducing the complexity of the congruent zeta function.

In the context of the string theoretic modularity problem the
L-function associated to the full cohomology group of a variety is
a physically useful object only in the case of elliptic curves. In
higher dimensions these functions are in general not modular. This
motivates the consideration of a factorization of L$-$functions,
and to ask whether modular forms emerge from the emerging pieces,
and if so, whether these modular forms admit a Kac-Moody theoretic
interpretation.

It was shown in refs. \cite{s06, s07b}, in the context of
establishing string theoretic modularity of certain Calabi-Yau
varieties, that a useful way to factorize the cohomological
L-functions of Calabi-Yau varieties is to consider subspaces of
the cohomology defined by Galois orbits of the holomorphic
$n-$forms $\Omega$, $n={\rm dim}_{\mathbb C}X$, where the Galois
groups are defined in an inherent way by the arithmetic of the
variety. The resulting L-functions have rational coefficients, and
therefore can in principle admit a string theoretic interpretation
along the lines discussed in \cite{su02, ls04, s05, s06, s07,
s07b}. In this paper we generalize this strategy to the case of
Fano varieties of special type by considering the orbit motive of
a Galois representation determined by the generator of the
cohomology group $H^{n-(Q-1),Q-1}(X)$ with $Q>1$. We continue to
denote the generator of
 this group by $\Omega$, and denote the motive by $M_{\Omega}(X)$. The goal
 is to compute the L-series associated to this motive, and
 to find a string theoretic interpretation of this function, following
 the strategy used in \cite{s06, s07b}. More precisely, we need to
 consider motives with Tate twists, introduced by Tate in the context of
 the Tate conjectures, which assert that certain subgroups of
 the $i^{\rm th}$ twists $H^{2i}(X)(i)$ of the cohomology groups
 $H^{2i}(X)$ are generated by classes of algebraic cycles \cite{t65,t94}
 (see also \cite{jps87}).
  In the present context we consider for any special Fano variety of
 charge $Q$ the twisted omega motive $M_{\Omega}(X)(Q-1)$ and its L-function,
 which we will denote by
 $L_{\Omega}^{Q-1}(X,s)=L(M_{\Omega}(X)(Q-1),s)$. We will
 call $M_{\Omega}(X)(Q-1)$ the $\Omega-$motive with a Tate twist of
 charge $(Q-1)$, or simply the omega motive of charge $(Q-1)$

In analogy to the general considerations in \cite{s06,s07b} for
the class of Calabi-Yau varieties in arbitrary dimensions, we can
 ask the following for the class of Fano varieties of special
 type.

\noindent
 {\bf Question 1.}
   Is the L-function $L_{\Omega}^{Q-1}(X,s)$ of the
   $\Omega-$motive $M_{\Omega}(X)(Q-1)$ of charge $(Q-1)$
   of a special Fano variety of
   charge $Q$ always modular?

 \noindent
 {\bf Question 2.} Can modular inverse Mellin transforms of
 $L_{\Omega}^{Q-1}(X,s)$ be expressed in terms of Hecke indefinite
 modular forms associated to Kac-Moody algebras?

It will become clear below that the answer to this question is
affirmative, at least for certain examples. In such cases the
L-function can be viewed as a map that takes motives and turns
them into conformal field theoretic objects. More challenging, in
part because of the lack of definition of some of the terms,
therefore would be the following picture.  \hfill \break
 {\bf Conjecture:}  \hfill \break
 The L-function is a functor from the additive category of motives
 of special Fano type to the multiplicative category of
 $N=2$ supersymmetric conformal field theories.

The basic outline described below is a generalization of the
method based on Jacobi sums of weighted projective spaces
described in \cite{s06} in the context of K3 surface modularity,
and in \cite{s07b} for higher dimensional Calabi-Yau varieties.

\vskip .3truein

\section{Modularity of Hecke L-series}

The modularity of the L-series determined in this paper follows
from the fact that they can be interpreted in terms of Hecke
L-series associated to Gr\"o\ss encharacters, defined by Jacobi
sums. A Gr\"o\ss encharacter can be associated to any number
field. Hecke's modularity discussion \cite{h37b} has been extended
by Shimura \cite{s71b} and Ribet \cite{r77}.

Let $K$ be a number field and
  $\sigma:~K ~\longrightarrow ~{\mathbb C}$ denote an embedding.

{\bf Definition.}~ {\it A Gr\"o\ss encharacter is a homomorphism
$\psi: I_{{\mathfrak m}}(K) \rightarrow {\mathbb C}^{\times}$ of
the
 fractional ideals of $K$ prime to the congruence integral ideal ${\mathfrak m}$
 such that
   $\psi((z)) = \sigma(z)^{w-1}$, $\forall K^{\times} \ni z
   \equiv 1({\rm mod}^{\times}{\mathfrak m})$.}

In the present case the cyclotomic Jacobi sums determined by the
finite field Jacobi sums computed above arise from imaginary
quadratic fields, in which case the structure of these characters
simplifies. Consider an imaginary quadratic field $K$ with
discriminant $-D$ and denote by $\varphi$ the Dirichlet character
associated to $K$, viewed as a mod $D$ character. Define a second
Dirichlet character $\lambda$ defined mod ${\rm N}{\mathfrak m}$
by setting
  \begin{equation}
   \lambda(a) = \psi((a))/\sigma(a)^{w-1},\phantom{gap}a\in {\mathbb Z}.
  \end{equation}
 The product of $\varphi$ and $\lambda$ then determines the
 nebentypus character of the modular form determined by Hecke character
 $\psi$.

Denote by ${\rm N}{\mathfrak p}$ the norm of prime ideal
${\mathfrak p}$, and by
 \begin{equation}
  L(\psi,s) = \prod_{{\mathfrak p} \in {\rm Spec}~{\mathcal O}_K}
  \left(1- \psi({\mathfrak p}){\rm N}{\mathfrak
  p}^{-s}\right)^{-1}
  \end{equation}
 the Hecke L-series associated to the character $\psi$.
 The modularity of the corresponding $q-$series
$f(\psi,q)=\sum_n a_nq^n$ associated to the L-series
  via the Mellin transform is characterized by the following result
  of Hecke.

  {\bf Theorem 4.}~
  {\it The power series $f(\psi,q)$ is a cusp form of weight $w$ and
   character $\epsilon = \varphi \lambda$ on $\Gamma_0(D{\rm N} {\mathfrak m})$.
   If $p {\relax{~|\kern-.35em /~}} D{\rm N}{\mathfrak m}$ then
     $f|_{T_p} = a_p f$, where $T_p$ are the Hecke operators.}

\section{Special Fano varieties with critical dimension three}

The focus in this paper will be on Fano varieties of
Brieskorn-Pham type, hence the computation of their L-functions
can proceed via their Jacobi sums \cite{w49,w52}.

{\bf Theorem 5.}  {\it For a smooth weighted projective
hypersurface with degree vector ${\underline n}=(n_0,...,n_s)$}
 \begin{equation}
 X^{\underline n} =
  \{b_0z_0^{n_0}+b_1 z_1^{n_1} + \cdots + b_s z_s^{n_s} =0 \}
  \subset {\mathbb P}_{(k_0,k_1,\dots,k_s)},
   \end{equation}
  {\it defined over the finite field ${\mathbb F}_q$, for $q=p^r, r\in
  {\mathbb N}$,
   set $d_i=(n_i,q-1)$  and denote by ${\mathcal A}_s^{q,{\underline n}}$
   the set of rational vectors
  $\alpha =(\alpha_0,\alpha_1,\dots,\alpha_s)$ given by}
  \begin{equation}
  {\mathcal A}_s^{q,{\underline n}}
  = \left\{\alpha \in {\mathbb Q}^{s+1} ~{\Big |}~
  0<\alpha_i<1, ~d_i\alpha_i =0~
  {\rm mod}~1,~\sum_{i=0}^s \alpha_i=0~({\rm mod}~1) \right\}.
  \end{equation}
  {\it For each $(s+1)-$tuple $\alpha$ define the Jacobi sum }
 \begin{equation}
 j_q(\alpha) = \frac{1}{q-1}
  \sum_{\stackrel{u_i \in {\mathbb F}_q}{u_0+u_1+\cdots +u_s=0}}
   \chi_{\alpha_0}(u_0) \chi_{\alpha_1}(u_1) \cdots \chi_{\alpha_s}(u_s),
  \end{equation}
  {\it where $\chi_{\alpha_i}(u_i) = e^{2\pi i \alpha_i m_i}$ with integers $m_i$
 determined via $u_i = g^{m_i}$, where $g\in {\mathbb F}_q$ is a
generator. Then the cardinality of $X^{\underline n}/{\mathbb
F}_q$ is given by}
 \begin{equation}
  N_{q}(X^{\underline n}) = 1+q + \cdots + q^{s-1}+
  \sum_{\alpha \in {\mathcal A}_s^{q,{\underline n}}} j_q(\alpha).
 \end{equation}

With these Jacobi sums one defines for primes $p$ such that
$d|(p^f-1)$ for some $f\in {\mathbb N}$ the polynomials
 \begin{equation}
 {\mathcal P}^{s-1}_p(t)
  = \prod_{\alpha \in {\mathcal A}_s^n}
    \left(1-(-1)^{s-1}j_{p^f}(\alpha)\prod_i{\bar \chi}_{\alpha_i}(b_i) t
   \right)^{1/f}
 \end{equation}
  and the associated L-function is given as a product over all
  good primes $p$
 \begin{equation}
 L^{(j)}(X,s) =
  \prod_p ({\mathcal P}_p^j(p^{-s}))^{-1}.
 \end{equation}
 The completion of the L-function at the bad primes will not be of
 relevance in the following.

The two  varieties considered here in the context of $D_{\rm
crit}$, i.e. Gepner models with central charge $c=9$, are
odd-dimensional, hence the rational form of the congruent zeta
function for a projective hypersurface of odd$-$dimension $s$
takes the form
 \begin{equation}
  Z(X_{s-1}/{\mathbb F}_p,t) = \frac{{\mathcal P}_p^{s-1}(t)}{(1-t)(1-pt)\cdots
  (1-p^{s-1}t)}.
 \end{equation}
We also consider a Fano variety corresponding to a K3 surface,
i.e. a Gepner model with central charge $c=6$. In this case the
smooth Fano variety is of complex dimension four and the zeta
function is of the form
 \begin{equation}
  Z(X^3_4/{\mathbb F}_p,t) = \frac{1}{(1-t)(1-pt){\mathcal P}_p^4(t)
  (1-p^3t)(1-p^4t)}.
 \end{equation}

 In all these examples there is therefore only one interesting
 cohomological L-function. It turns out that these L-functions are
 not themselves of interest and that it is more useful to consider the
  L-function of the $\Omega-$motive. In all three examples
 the order of the Galois group is two, hence the motive is of
 dimension two. Denoting the Jacobi sum associated to the $\Omega-$form
 by $j_p(\alpha_{\Omega})$ the polynomial of the $\Omega-$motive takes the
 form
 \begin{equation}
  {\mathcal P}_p^{\Omega}(t) = \left(1-j_{p^f}(\alpha_{\Omega})t^f\right)^{1/f}
                  \left(1-j_{p^f}({\bar
                  \alpha}_{\Omega})t^f\right)^{1/f},
 \end{equation}
 where ${\bar \alpha}_{\Omega}:={\bf 1} - \alpha_{\Omega}$ is the dual
 vector and ${\bf 1}$ is the unit vector.
  For primes $p$ such that $d|(p-1)$ this simplifies and the local L-function
 factors of the twisted omega motive of charge $(Q-1)$ take the form
 $L_{\Omega,p}^{Q-1}(X,s) = ({\mathcal P}_p^{\Omega,Q-1}(p^{-s}))^{-1}$
 with
 $${\mathcal P}_p^{\Omega,Q-1}(t) = 1 + \beta(p) p^{-(Q-1)}t +
 p^{n-2(Q-1)}t^2,$$
  where $\beta(p) = j_p(\alpha_{\Omega}) + j_p({\bar \alpha}_{\Omega})$.

\subsection{The cubic Fermat sevenfold $X_7^3$}

 \subsubsection{The motivic L-function}
 For the cubic Fermat sevenfold the Hodge decomposition takes the form
 $h^{i,i}(X_7^3) = 1, i=0,1,2,3$ with all other $h^{i,j}(X_7^3)=0$
 for $i,j<4, i\neq j$. The intermediate cohomology
 decomposes as in (\ref{cubic-sevenfold-hodges}), i.e. we have
 $Q=3$ and $D_{\rm crit}=3$. The zeta function therefore
 has the form $Z(X_7^3/{\mathbb F}_p,t) =
 {\mathcal P}_p^7(t)/\prod_{i=0}^7 (1-p^it)$,
  leading to the L-function
  \begin{equation}
   L(X_7^3,s) := \prod_p 1/{\mathcal P}_p^7(p^{-s}).
 \end{equation}

The cyclotomic field defined by this manifold is ${\mathbb
Q}(\mu_3)$, hence we have ${\rm Gal}({\mathbb Q}(\mu_3)/{\mathbb
Q}) = ({\mathbb Z}/3{\mathbb Z})^{\times}$, and the L-series of
the $\Omega-$motive of $X_7^3$ can therefore be computed from the
orbit of the Jacobi sums $j_p(\alpha_{\Omega})$ of length two,
with
 $\alpha_{\Omega} = \frac{1}{3}\left(1,1,1,1,1,1,1,1,1\right)$.
 The L-series of the twisted $\Omega-$motive of charge two of
 $X_7^3$ is defined as
 $L_{\Omega}^2(X_7^3,s) = \prod_p 1/{\mathcal P}_p^{\Omega,2}(p^{-s})$
  with ${\mathcal P}_p^{\Omega,2}(t)= 1 + a_p t + p^3t^2$
 for $n|(p-1)$, $\beta(p) = (j_p(\alpha_{\Omega}) +
       {\bar j}_p(\alpha_{\Omega}))$ and
 $a_p = \beta(p)/p^2$ because the Tate twist is of charge two.
 The computation of these sums for the first few nontrivial primes
 is sufficient to
 determine the L-function of the Tate twisted $\Omega-$motive
 \begin{equation}
 L_{\Omega}^2(X_7^3,s) \doteq 1 + \frac{20}{7^s} - \frac{70}{13^s}
  + \frac{56}{19^s} + \frac{308}{31^s} + \cdots
 \end{equation}
 and the associated $q-$series
 \begin{equation}
   f_{\Omega}^2(X_7^3,q) \doteq q + 20q^7 - 70q^{13} + 56q^{19} + 308q^{31}  + \cdots
 \end{equation}
 completely.

It will become clear below that this $q-$series describes a
modular cusp form of weight four and modular level nine with
respect to the Hecke congruence group.

It was suggested in \cite{s92,cdp93,s94} that the cubic Fermat
sevenfold can be interpreted as a geometric model associated to
the tensor model $1^{\otimes 9}$ of nine copies of $N=2$
superconformal models \cite{g88}. If this is indeed the case, then
we would expect, in analogy to the results obtained in \cite{su02,
ls04, s05, s06, s07, s07b}, that we can express the $q-$series of
the $\Omega-$motive in terms of the Hecke indefinite modular forms
$\Theta^k_{\ell,m}(\tau)$ associated to the underlying conformal
field theory model. At the affine level $k=1$ there is a single
such form, hence this is a make or break situation which is very
constrained, much like the first example established in
 this framework \cite{su02}. It turns out that the geometric modular form of the
$\Omega-$motive of the cubic sevenfold does admit a string
theoretic interpretation. There are several ways this can be
described. The most direct is that the modular form
$f_{\Omega}^2(X_7^3,q)$ can be written as
 \begin{equation}
 f_{\Omega}^2(X_7^3,q) = \Theta^1_{1,1}(q^3)^4 ~\in ~S_4(\Gamma_0(9)).
 \end{equation}
 This modular form has the expected weight
 of a Calabi-Yau threefold, consistent with the central charge of the
 model. The absence of a twist character in this result
 is consistent with the fact that the field of quantum dimensions
 of the underlying field theory is trivial \cite{s01}.

\subsubsection{CM interpretation of the L-series and modularity}

The motivic modular form $f_{\Omega}^2(X_7^3,q)$ of the cubic
sevenfold is a form with complex multiplication by the Eisenstein
field ${\mathbb Q}(\sqrt{-3})$, and therefore can be described by
a Hecke character associated to this field. Define a character
$\psi_{27}$ by the following congruence condition
 \begin{equation}
  \psi_{27}({\mathfrak p}) = \alpha_{\mathfrak p},
  \phantom{gap} \alpha_{\mathfrak p}
  \equiv 1({\rm mod}~3),
 \label{eischar1}\end{equation}
 where the generator of the prime ideal
 ${\mathfrak p}=(\alpha_{\mathfrak p})$ is
 determined uniquely by the congruence constraint.
  As described in the previous section, associated to any Hecke character
  $\psi$ is a Hecke series $L(\psi,s)$, and we find
  \begin{equation}
   L_{\Omega}^2(X_7^3,s) = L(\psi_{27}^3,s).
  \end{equation}
  It follows from Hecke's theorem that the $q-$series
  of this L-function is  modular, hence the data provided above is sufficient to prove
  that the Mellin transform of the L-function
  $L_{\Omega}^2(X_7^3,s)$ is a modular
  form, and therefore the Tate twisted $\Omega-$motive of the cubic sevenfold
  is modular.

This result shows that the basic arithmetic building block of the
cubic sevenfold is the elliptic cubic Fermat curve
 \begin{equation}
  E^3 = \left\{(z_0:z_1:z_2)\in {\mathbb P}_2~|~
  z_0^3+z_1^3+z_2^3=0\right\}.
 \end{equation}
 The Hasse-Weil L-series of $E^3$ agrees with the L-series of the
 algebraic Hecke character $\psi_{27}$
  \begin{equation}
   L_{\rm HW}(E^3,s) = L(\psi_{27},s).
  \end{equation}
  The Mellin transform of $L_{\rm HW}(E^3,s)$ has weight 2 and
  level 27, and it factors into string theoretically induced modular forms
  as \cite{su02}
  \begin{equation}
   f(E^3,q) = \Theta^1_{1,1}(q^3)\Theta^1_{1,1}(q^9).
 \end{equation}
 This proves the part of Theorem 1 concerned with the cubic
 Fermat sevenfold.

 The relation between twisted $\Omega-$motivic L-series of $X_7^3$ and the
 L-series of the elliptic curve $E^3$ that follows from these
 considerations can be formulated in a more direct way.
 Consider the expansions $L_{\Omega}^2(X_7^3,s) = \sum_n a_n n^{-s}$
 for the sevenfold,
 and $L_{\rm HW}(E^3,s) = \sum_n b_nn^{-s}$ for the elliptic curve.
 The Hecke character interpretation means that
  $a_p = \alpha_{\mathfrak p}^3 + {\bar \alpha}_{\mathfrak p}^3$ and
 $b_p = \alpha_{\mathfrak p} + {\bar \alpha}_{\mathfrak p}$.
 Therefore
  \begin{equation}
   a_p = b_p^3 -3pb_p.
  \end{equation}

\subsubsection{$X_7^3$ as the mirror of a rigid Calabi-Yau
variety}

The fact that the CM character of the elliptic curve $E^3\subset
{\mathbb P}_2$ is the building block of the L-function of the
cubic sevenfold suggests that its $\Omega-$motivic L-function
should be identical to that of the resolution of a quotient of the
triple product $(E^3)^3$ which produces the correct mirror
cohomology of the critical part of the spectrum of the sevenfold.
Such a quotient variety can be obtained by considering the
 ${\mathbb Z}_3\times {\mathbb Z}_3$ action defined on $(E^3)^3$ as
 \begin{equation}
  {\mathbb Z}_3\times {\mathbb Z}_3 \ni g_1\times g_2:~~\left[
  \begin{matrix}
    1 &0 &0 &0 &0 &0 &2 &0 &0\cr
    0 &0 &0 &1 &0 &0 &2 &0 &0\cr
  \end{matrix}
  \right]
 \end{equation}
 with obvious notation for the generators $g_1,g_2$. The
 resolution $X_3^3$ of the quotient variety
  $(E^3)^3/{\mathbb Z}_3\times {\mathbb Z}_3$ is rigid and has
  $h^{1,1}(X_3^3)=84$, as needed for the mirror cohomology of the
  relevant part of the cohomology of the sevenfold. Higher
  dimensional varieties defined by the resolution of quotients of
  the type $E^n/G$, for elliptic curves and certain finite groups $G$ were
  considered by Cynk and Hulek \cite{ch05}.

  {\bf Theorem 6.}~
  {\it Let $n\in {\mathbb N}$ be an odd integer,
  $E$ be the elliptic curve with an automorphism
  of order 3, and denote by ${\bar X}_n$ the quotient of $E^n$ by the action
  of the group }
  $$
  \left\{(\xi^{a_1}\times \cdots \times \xi^{a_n}) \in
  {\rm End}(E^n)~|~
   \sum_{i=1}^n a_i = 0({\rm mod} ~3)\right\}.
  $$
  {\it Then ${\bar X}_n$ has a smooth model $X_n$, which is a
  Calabi-Yau manifold such that ${\rm dim}~H^n(X_n) = 2$
  and $L(X_n,s) \stackrel{\cdot}{=} L(g_{n+1},s)$, where $g_{n+1}$
  is the weight $n+1$ cusp form with complex
  multiplication by ${\mathbb Q}(\sqrt{-3})$, associated to the
 $n^{\rm th}$ power of the Gr\"o\ss encharacter of $E$.}

 It follows by combining this result with those of Theorem 1 that the
  L-function of the $\Omega-$motive of the cubic sevenfold,
  with a Tate twist of charge two, is identical to the L-function of
  the intermediate cohomology of the variety $X_3^3$. This completes
  the proof of the part of Theorem 2
  concerned with the cubic Fermat sevenfold.

\subsubsection{Relatives and correspondences}

The agreement between the motivic L-function of $X_7^3$ and the
L-function of the intermediate cohomology of $X_3^3$ identifies
these two varieties as "relatives" in the sense of ref.
\cite{hsvgvs01}. A variety $X_1$ with a 2-dimensional Galois
representation $\rho_1$ is defined to be a relative of a variety
$X_2$ if $\rho_1$ occurs in the cohomology of $X_2$. This notion
identifies $X_3^3$ as a relative of $X_7^3$. Other relatives of
Calabi-Yau type for $X_3^3$ have been considered by van Geemen in
\cite{wvg90}, Saito-Yui01 \cite{sy01}, Verrill \cite{v03}, and Yui
\cite{y03}.

A stronger relation is predicted by the Tate conjecture, which in
the present context can be formulated as a statement that an
isomorphism between two 2-dimensional Galois representations which
occur in the \'etale cohomology of two varieties defined over
${\mathbb Q}$ should be induced by a correspondence between the
two varieties defined over ${\mathbb Q}$, i.e. there should be an
algebraic cycle on the product of the two varieties \cite{y03}.
The Tate conjecture therefore implies that there should exist a
correspondence between $X_3^3$ and $X_7^3$. An explicit
construction of such a correspondence should provide better
insight into the precise relation between the pair of manifolds
$(X_7^3, X_3^3)$.

The fact that physically the L-function is determined by the
modular form $\Theta^1_{1,1}(\tau)$ of the underlying conformal
field theory on the worldsheet leads to the conclusion that from a
fourdimensional point of view these two varieties contain motivic
sectors that are indistinguishable.

\subsection{The quartic Brieskorn-Pham fivefold $X_5^4$}

 \subsubsection{The L-series of the $\Omega-$motive}

 The quartic fivefold $X_5^4$ defined in (\ref{fano2})
  has a Hodge decomposition of the form
 $h^{i,i}(X_5^4) = 1, i=0,1,2$, with all other
 $H^{i,j}(X_5^4)=0$ for $i,j<5$ and $i\neq j$.
 The decomposition of the intermediate cohomology group is
 given by (\ref{quartic-fivefold-hodges}),
  i.e. $Q=2$ and $D_{\rm crit}=3$.

 The congruence zeta function therefore takes the form
 $Z(X_5^4/{\mathbb F}_p,t) =
  {\mathcal P}_p^5(t)/\prod_{i=0}^5(1-p^it)$,
 leading to the L-function
   $L(X_5^4,s) = \prod_p ({\mathcal P}_p^5(p^{-s}))^{-1}$.
  We are interested in the factor describing the twisted
  $\Omega-$motive of charge one, which is determined by the Jacobi sums
  $j_p(\alpha_{\Omega})$ with
  $ \alpha_{\Omega} = \frac{1}{4}(1,1,1,1,1,1,2)$.
 The computation of the first few Jacobi sums for low primes leads
 to the $q-$expansion of the L-function of the twisted $\Omega-$motive
of the quartic fivefold given by
 \begin{equation}
    f_{\Omega}^1(X_5^4,q) \doteq q - 22q^5 + 18q^{13} -  94q^{17} + 130q^{29} + \cdots
 \end{equation}
 This is a modular cusp form of weight four and level 64,
  $f_{\Omega}^1(X_5^4,q) \in S_4(\Gamma_0(64))$.

 \subsubsection{CM interpretation of the L-series and modularity}

 The motivic modular form $f_{\Omega}^1(X_5^4,q)$
 has CM by the Gauss field ${\mathbb Q}(\sqrt{-1})$, and therefore can be
described by a Hecke character associated to this field as well.
Consider the twisted character $\psi_{64}=\psi_{32} \chi_2$ where
$\psi_{32}$ is defined by
 \begin{equation}
  \psi_{32}({\mathfrak p}) = \alpha_{\mathfrak p},\phantom{gap}
     \alpha_{\mathfrak p} \equiv 1({\rm mod}~(2+2i)),
 \end{equation}
 where the generator of the prime ideal ${\mathfrak p}=(\alpha_{\mathfrak p})$ is
 determined uniquely by the constraint.

  The Hecke series associated to $\psi_{64}$ then is again the building block
  of the geometric L-series
  \begin{equation}
   L_{\Omega}^1(X_5^4,s) = L(\psi_{64}^3,s)
  \end{equation}
 hence modularity of this L-function follows. Similar to the cubic
 sevenfold this result shows that the quartic fivefold has an
 elliptic building block, in this case the Brieskorn-Pham type
curve
 \begin{equation}
  E^4 = \left\{(z_0:z_1:z_2)\in {\mathbb P}_{(1,1,2)}~{\Big |}~
   z_0^4+z_1^4+z_2^2=0\right\}.
  \end{equation}

 It can be checked that the Hasse-Weil L-series of $E^4$ is
 determined by the character $\psi_{64}$
 \begin{equation}
  L_{\rm HW}(E^4,s) = L(\psi_{64},s)
 \end{equation}
 and it was shown in \cite{ls04} that the Mellin transform of
 $L_{\rm HW}(E^4,s)$ factors into string theoretic Hecke
 indefinite modular forms as
 \begin{equation}
  f(E^4,q) = \Theta^2_{1,1}(q^4)^2 \otimes \chi_2.
  \end{equation}
 This completes the proof of Theorem 1.

 \subsubsection{A rigid Calabi-Yau mirror}

The appearance of the quartic elliptic curve $E^4$ again suggests
to consider a resolved quotient of the triple product $(E^4)^3$,
defined by some discrete symmetry group action such that the
resolved manifold has a Hodge decomposition that is mirror to the
critical part of the spectrum of the quartic fivefold. It turns
out that the group that produces such a resolution is given by
 ${\mathbb Z}_4\times {\mathbb Z}_4$, with the action
 \begin{equation}
  {\mathbb Z}_4\times {\mathbb Z}_4 \ni g_1\times g_2:~~\left[
  \begin{matrix}
    1 &0 &0 &0 &0 &0 &3 &0 &0\cr
    0 &0 &0 &1 &0 &0 &3 &0 &0\cr
  \end{matrix}
  \right].
 \end{equation}
 Denote by $X_3^4$ the resolution of the quotient variety
  $(E^4)^3/{\mathbb Z}_4\times {\mathbb Z}_4$. This is a rigid
  Calabi-Yau variety which has
  $h^{1,1}(X_3^4)=90$, as needed for the mirror cohomology of the
  relevant part of the cohomology of the fivefold. Again we can
  apply  result of Cynk and Hulek \cite{ch05}
  to the mirror just constructed.

 {\bf Theorem 7.}~
  {\it Let $n\in {\mathbb N}$ be an odd integer,
  $E$ be the elliptic curve with an automorphism
  of order 4, and denote by ${\bar X}_n$ the quotient of $E^n$ by the action
  of the group }
  $$
  \left\{(\xi^{a_1}\times \cdots \times \xi^{a_n}) \in
  {\rm End}(E^n)~|~
   \sum_{i=1}^n a_i = 0({\rm mod} ~4)\right\}.
  $$
  {\it Then ${\bar X}_n$ has a smooth model $X_n$, which is a
  Calabi-Yau manifold such that ${\rm dim}~H^n(X_n) = 2$
  and $L(X_n,s) \stackrel{\cdot}{=} L(g_{n+1},s)$, where $g_{n+1}$
  is the weight $n+1$ cusp form with complex
  multiplication by ${\mathbb Q}(\sqrt{-1})$, associated to the
 $n^{\rm th}$ power of the Gr\"o\ss encharakter of $E$.}

  It follows from Theorems 1 and 7 that the L-functions of
  the quartic fivefold and the variety $X_3^4$
  are identical. This completes the proof of Theorem 2.

\subsubsection{Relatives and correspondences}

As in the case of the mirror pair $(X_7^3,X_3^3)$, the pair of
varieties $(X_5^4, X_3^4)$ form a pair of relatives in that the
Galois representation of $X_3^4$ is contained in the cohomology of
the quartic fivefold. Again, the Tate conjecture implies that
there is a correspondence between $X_5^4$ and $X_3^4$. The fact
that the L-function has a physical interpretation in terms of the
worldsheet modular form $\Theta^2_{1,1}(\tau)$ leads again to the
conclusion that from a low energy perspective the corresponding
motives of these two varieties are indistinguishable.

\vskip .2truein

\section{A special Fano variety of critical dimension two}

\subsection{Arithmetic of the cubic Fermat fourfold $X_4^3$}

The relation between special Fano varieties and critical string
models manifolds exists in all dimensions. An interesting
lower-dimensional case is that of K3 surfaces. Among the Gepner
models there exists an example that is the direct analog of the
cubic 7$-$fold, defined by the cubic fourfold
 (\ref{cubic-fourfold}).
  The Hodge diamond structure of this variety is
  $h^{i,i}(X_4^3) = 1,~i=0,1$, with all other $H^{i,j}(X_4^3)=0$
  for $i,j<4$ and $i\neq j$, while the intermediate
  cohomology group decomposes as
 described above in eq. (\ref{cubic-fourfold-hodges}),
 i.e. $Q=2$ and $D_{\rm crit}=2$.

 The middle dimensional cohomology therefore is that of the K3
 surface, modulo the sector inherited from the ambient space ${\mathbb P}_5$.
 The zeta function of this cubic fourfold therefore takes the form
 \begin{equation}
  Z(X_4^3,t) = \frac{1}{(1-t)(1-pt){\mathcal
  P}_4^{(p)}(t)(1-p^3t)(1-p^4t)},
 \end{equation}
  and it is natural to ask whether the $\Omega-$motive of this
  variety leads to the L-function of an $\Omega-$motive of K3-type.

 The Galois group associated to this variety is again
 $({\mathbb Z}/3{\mathbb Z})^{\times}$ and the computation of the
 Tate twisted $\Omega-$motive obtained from the Jacobi sums
 $j_p(\alpha_{\Omega})$ with $\alpha_{\Omega} = \frac{1}{3}(1,1,1,1)$
 leads to $q-$expansion
 \begin{equation}
  f_{\Omega}^1(X_4^3,q) \doteq q - 13q^7 - q^{13} + 11q^{19} - 46q^{31} + \cdots
 \label{formofcubic4fold}
 \end{equation}
 This series describes a weight 3 modular form of level 27 which
 can be described in closed form as
  \begin{equation}
   f_{\Omega}^1(X_4^3,q) =  \vartheta(q^3)\Theta^1_{1,1}(q^3)\Theta^1_{1,1}(q^9)
   \in S_3(\Gamma_0(27)),
  \end{equation}
  where $\vartheta(q)$ is the Eisenstein series
  \begin{equation}
    \vartheta(q) = \sum_{z\in {\mathcal O}_K} q^{{\rm N} z}
  \end{equation}
   associated to the field $K={\mathbb Q}(\sqrt{-3})$.

 It was shown in \cite{s06} that (\ref{formofcubic4fold})
  is precisely the form which
   describes the Mellin transform of the L-series of the
  $\Omega-$motive of the degree six K3 surface
  \begin{equation}
  S^6 = \left\{(z_0:\cdots z_3)\in {\mathbb P}_{(1,1,1,3)}
   ~{\Big |} ~z_0^6 + z_1^6 + z_2^6 + z_3^2 = 0\right\}.
 \end{equation}
  Hence we have
  \begin{equation}
    f_{\Omega}^1(X_4^3,q) = f_{\Omega}(S^6,q).
 \end{equation}
 It is therefore natural to ask whether these two varieties can be
 related in some explicit way. This can be analyzed via the twist map.

 \subsection{Twist map construction}

 In the notation of \cite{hs99,hs95} consider the map
 \begin{equation}
 \Phi: ~{\mathbb P}_{(w_0,...,w_m)} \times {\mathbb P}_{(v_0,...,v_n)}
  \longrightarrow {\mathbb P}_{(v_0w_1,...,v_0w_m,w_0v_1,...,w_0v_n)}
 \end{equation}
 defined as
  \begin{eqnarray}
  ((x_0,...,x_m),(y_0,...,y_n) &\mapsto &
   (y_0^{w_1/w_0}x_1,...,y_0^{w_m/w_0}x_m,
   x_0^{v_1/v_0}y_1,...,x_0^{v_n/v_0}y_n) \nonumber \\
   ~~~& & =:(z_1,...,z_m,t_1,...,t_n).
   \end{eqnarray}
 This map restricts on the subvarieties
  \begin{eqnarray}
   X_1 &=& \{x_0^{\ell} + p(x_i) =0\} \subset {\mathbb P}_{(w_0,...,w_m)}
   \nonumber \\
   X_2 &=& \{y_0^{\ell} + q(y_j) =0\} \subset {\mathbb P}_{(v_0,...,v_n)}
  \end{eqnarray}
  to the hypersurface
  \begin{equation}
   X = \{p(z_i) - q(t_j) =0\} \subset
    {\mathbb P}_{(v_0w_1,...,v_0w_m,w_0v_1,...,w_0v_n)}
  \end{equation}
  as a finite map.
  The degrees of the hypersurfaces $X_i$ are given by
   \begin{equation}
    {\rm deg}~X_1 = w_0\ell,~~~~~~{\rm deg}~X_2 = v_0\ell,
  \end{equation}
   leading to the degree ${\rm deg}~X = v_0w_0\ell$.

Applying this construction to the hypersurfaces
 \begin{equation}
  S^6_{\pm} =
   \left\{(z_0:\cdots z_3)\in {\mathbb P}_{(1,1,1,3)}
   ~{\Big |} ~z_0^6 + z_1^6 + z_2^6 \pm z_3^2 = 0\right\}
 \end{equation}
  leads to the cubic fourfold embedded in ${\mathbb P}_5$, explaining
  geometrically the L-function result.

\subsection{CM interpretation of the L-series and modularity}

The modular form $f_{\Omega}^1(X_4^3,q)$ of the twisted omega
motive of the cubic fourfold has complex multiplication by the
Eisenstein field ${\mathbb Q}(\sqrt{-3})$ with the same underlying
basic character $\psi_{27}$ defined in (\ref{eischar1}) as the
cubic sevenfold, except that in the present case we have to
consider the square of the character
 \begin{equation}
  L_{\Omega}^1(X_4^3,s) = L(\psi_{27}^2,s).
 \end{equation}
 Hence the $\Omega-$motive of the cubic fourfold is modular. The
 basic arithmetic building block is again the elliptic curve
 $E^3$, described in string theoretic terms in \S6.

\end{document}